\documentclass[12pt]{article}
\usepackage{authblk}
\usepackage{hyperref}
\usepackage{amssymb}
\usepackage{amsthm}
\usepackage{amsmath}
\usepackage{multirow}
\usepackage{float}
\usepackage{graphicx}
\usepackage{listings}
\usepackage{MnSymbol}
\usepackage{nicematrix}
\usepackage{titlesec}
\usepackage[a4paper,margin=1in]{geometry}
\usepackage{forest}
\usepackage{xcolor}
\usepackage{tikz}
\usepackage{relsize}

\usetikzlibrary{trees, positioning, shapes.geometric}
\usetikzlibrary{trees}
\newtheorem{theorem}{Theorem}[section]

\newtheorem{example}[theorem]{Example}

\newtheorem{remark}[theorem]{Remark}

\begin{document}
\title{\Large\bfseries Invariant characterization of scalar third-order ODEs that admit the maximal contact symmetry Lie algebra}

\author[1]{Omar A. Abuloha\thanks{1227001@student.birzeit.edu}}
\author[1]{Marwan Aloqeili\thanks{maloqeili@birzeit.edu}}
\author[1]{Ahmad Y. Al-Dweik\thanks{aaldweik@birzeit.edu}}
\author[2]{F. M. Mahomed\thanks{Fazal.Mahomed@wits.ac.za}}
\author[3]{M. T. Mustafa\thanks{Corresponding Author: tahir.mustafa@qu.edu.qa}}
\affil[1]{Department of Mathematics, Birzeit University, Ramallah, Palestine}
\affil[2]{School of Computer Science and Applied Mathematics, 
University of the Witwatersrand, Johannesburg, Wits 2050,  South Africa}
\affil[3]{Department of Mathematics and Statistics, College of Arts and Sciences, Qatar University, 2713, Doha, Qatar}
\date{}
\maketitle
\begin{abstract}
The Cartan equivalence method is utilized to deduce an invariant
characterization of the scalar third-order ordinary differential
equation $u'''=f(x,u,u',u'')$ which admits the maximal
ten-dimensional contact symmetry Lie algebra. The method provides
auxiliary functions which can be used to efficiently determine the
contact transformation that does the reduction to the simplest
linear equation $\bar{u}'''=0$. Furthermore, ample examples are given to
illustrate our method.
\end{abstract}
\section{Introduction}
The idea of tangential transformations was initiated in the early
work of Lie \cite{Lie1} wherein he found a local contact
transformation which mapped straight lines into spheres in space.
Notwithstanding, the theoretical foundation of the theory of
contact transformations is given in Lie and Engel
\cite{Lie2,Lie3,Lie4,Lie5}.  Lie showed that for scalar
$n$th-order ordinary differential equations (ODEs),  $n = 3$, its
contact symmetry algebra is of finite dimension. He presented the
complete classification of finite-dimensional irreducible contact
Lie algebras in two complex variables \cite{Lie4}.  In the more
recent paper \cite{Wafo}, the authors classified $n$th-order ODEs,
$n\ge 3$, that admit nontrivial contact symmetry Lie algebras. Lie
and Scheffers \cite{Lie6} have shown that a third-order ODE admits
at most a ten-dimensional contact symmetry Lie algebra. They
proved that the symmetry algebra is ten-dimensional if and only if
the third-order ODE  is equivalent, up to a local contact
transformation, to the simplest third-order ODE (see also
\cite{Wafo}). The reader is also referred to
\cite{Svi1,Svi2,Ibr,Yum1} on these and related aspects.

Cartan \cite{car}, inter alia, provided a solution to the linearization problem for
ODEs using his now popular approach called the Cartan equivalence method
(see the recent contributions \cite{Gardner1989,Olver1995}).

Lie and Scheffers \cite{Lie6}, Yumaguzhin \cite{Yum1,Yum} as well as Wafo Soh et al. \cite{Wafo} studied the local classification of third-order linear ODEs, up to contact transformations. There are three canonical forms that occur for scalar linear
third-order ODEs. The maximal contact symmetry Lie algebra case
for such ODEs is of dimension ten and corresponds to the simplest equation $u'''=0$.

The Laguerre-Forsyth  canonical form for scalar
linear third-order ODEs is given by (see \cite{mah4})
\begin{equation}\label{e3}
u'''+a(x)u=0,
\end{equation}
where $a(x)$ is an arbitrary function of $x$. If $a(x) \not \equiv
0$, then it is known that (\ref{e3}) has a five- or
four-dimensional contact symmetry algebra and otherwise
ten-dimensional contact symmetry algebra.

Chern \cite{che} was the first to invoke the Cartan equivalence
method in order to solve the linearization problem for scalar
third-order ODEs via contact transformations. He deduced
conditions of equivalence to the equation (\ref{e3}) in the two
cases $a(x)\equiv0$ and $a(x)\equiv1$. Then Neut and Petitot
\cite{neu} also studied equivalence to (\ref{e3}) by means of
contact transformations but for arbitrary $a(x)$.  In the recent work,
Ibragimov and Meleshko \cite{ibr1} investigated the linearization
problem for third-order ODEs by utilising a direct approach via
both point and contact transformations. It is the case that
Ibragimov and Meleshko \cite{ibr1} first studied the second part
of the linearization problem for scalar third-order ODEs which is
that of constructing the relevant transformations to simpler form,
via both point and contact transformations. However, it should be
mentioned that the solution of the second part was given as a
solution of a non-linear system of PDEs in their investigation. In
the works
\cite{DweikTahirFazal1,DweikTahirFazal2,DweikTahirFazal3}, the
authors very recently applied a new framework of the Cartan
equivalence method  to solve the two parts of the linearization
problem for scalar third-order ODEs of the form
$u'''=f(x,u,u',u'')$ via point transformations. In these essential
works, the transformations can be obtained efficiently as a
solution of a system of linear or Riccati equations given in terms
of the introduced auxiliary functions.

We emphasise here that an invariant characterization of a third-order ODE
$u'''=f(x,u,u',u'')$ which admits the maximal ten contact symmetry
algebra was obtained in terms of the function $f$ in the following
theorem. We denote $u', u''$ by $p, q$, respectively in the
following and in what transpires in the sequel.
\begin{theorem}\cite{neu}
The third-order equation $ u''' = f(x,u,u',u'')$ is equivalent to the simplest
form $u''' = 0$ with ten contact symmetries under {\rm contact
transformations} if and only if the relative invariants
\begin{equation}\label{wer2}
\begin{array}{l}
 K_1  = f_{q,q,q,q}  \hfill \\
 K_2  =  {4\,f_q ^3  + 18\,f_q \left( {f_p  - \hat{D}_x f_q } \right) + 9\,\hat{D}_x^2 f_q    - 27\,\hat{D}_x f_p+ 54\,f_u  }\\
\end{array}
\end{equation}
both vanish identically, where $ \hat{D}_x = \frac{\partial } {{\partial
x}} + p\frac{\partial }{{\partial u}}+q\frac{\partial} {{\partial
p}}+f\frac{\partial} {{\partial q}}$ and  $K_2$ is the well-known
W$\ddot{\textrm{u}}$nschmann relative invariant
\cite{neu,Wunschmann}.
\end{theorem}
In this paper, we apply Cartan's equivalence method to get an invariant coframe for the equivalence problem of third order ODEs that admit maximal contact symmetry Lie algebra, which is achieved in the next section. Then, we provide in section 3 a systematic framework for constructing the contact transformations that reduce the third-order ODE with ten contact symmetries to its canonical form using the obtained invariant coframe. It is opportune to
 remark that the application of the Cartan equivalence method
in this framework is both relevant and new. The new framework of Cartan's equivalence method gives invariant coframe explicitly in
terms of auxiliary functions. The invariant coframe is utilized to
determine the contact transformations to the equivalent canonical
form. In section 4, some important examples are provided to explain our results.
\section{Application of Cartan’s Equivalence Method to Third-Order ODEs with Maximal Contact Symmetry}
\setcounter{equation}{0}
In this section, Cartan's equivalence method is applied to derive an invariant coframe for scalar third-order ODEs admitting the maximal ten-dimensional contact symmetry Lie algebra. The reduction proceeds through successive normalizations of the structure group until an invariant e-structure is obtained.

The reader is referred to \cite{Gardner1989,Olver1995} for the basic definitions, notation, and results needed in this section.
Let $\left(x, u, p=u^{\prime},q=u^{\prime\prime}\right) \in \mathbb{R}^4$ be local coordinates on the second-order jet space $\bf{J}^2$.  Throughout this paper, the 1-forms $\pi^{\prime\kappa},~\kappa=1, \ldots, 9$,   denote the modified Maurer-Cartan forms. Consider the following base coframe on $M=\bf{J}^2$
\begin{equation}\label{2.1}
\left(\begin{array}{l}
\omega^1\\\omega^2\\\omega^3\\\omega^4\\
\end{array}\right)
=\left(\begin{array}{c}
du-p dx\\
dp-q dx\\
dq-f dx\\
dx\\
\end{array}\right),
\end{equation}
In local coordinates, the equivalence  of
\begin{equation}\label{2.3}
\begin{array}{l}
u^{\prime \prime \prime}=f\left(x, u, u^{\prime},u^{\prime \prime} \right), \quad \bar{u}^{\prime \prime\prime}=\bar{f}\left(\bar{x}, \bar{u}, \bar{u}^{\prime},\bar{u}^{\prime \prime}\right),
\end{array}
\end{equation}
under a contact transformation
\begin{equation}\label{2.4}
\bar{x}=\varphi \left( x,u,p \right),~\bar{u} =\psi \left( x,u,p  \right), ~\bar{p}=\chi(x,u,p),\\
\end{equation}
satisfying $d\bar{u}-\bar{p}~d\bar{x}=\lambda(du-p~dx)$ for some function $\lambda(x,u,p)$ and having a nonvanishing Jacobian, can be expressed in terms of the base coframe (\ref{2.1}) by the following conditions
\begin{equation}\label{2.5}
\Phi^*\left(\begin{array}{c}
\bar{\omega}^1 \\
\bar{\omega}^2 \\
\bar{\omega}^3\\
\bar{\omega}^4
\end{array}\right)=\left(\begin{array}{cccc}
a_1 & 0 & 0&0 \\
a_2 & a_3 & 0 &0\\
a_4 & a_5 & a_6&0\\
a_7&a_8&0&a_9
\end{array}\right)\left(\begin{array}{c}
\omega^1 \\
\omega^2 \\
\omega^3\\
\omega^4
\end{array}\right),
\end{equation}
for functions $a_i(x,u,p,q),~ i=1,...,9$,  where $\Phi^*$ is the pullback arising from the first prolongation of the contact transformation (\ref{2.4}). Consequently, the corresponding structure group is a nine-dimensional Lie group
\begin{equation}\label{2.5}
G=\left\{
\left(\begin{array}{cccc}
a_1 & 0 & 0&0 \\
a_2 & a_3 & 0 &0\\
a_4 & a_5 & a_6&0\\
a_7&a_8&0&a_9
\end{array}\right)  \Bigg \vert a_1a_3a_6a_9\neq 0
\right\}.
\end{equation}
Now, let us define $\theta$ to be the lifted coframe
\begin{equation}\label{2.6}
\begin{array}{ll}
\left(\begin{array}{c}
\theta^1, ~\theta^2,~\theta^3,~\theta^4
\end{array}\right)^T=S \left(\begin{array}{c}
\omega^1,~
\omega^2,~
\omega^3,~
\omega^4
\end{array}\right)^T,
\end{array}
\end{equation}
where $S\in G$. 
After absorption, the first structure equation for (\ref{2.6}) is
\begin{equation}\label{2.7}
\begin{array}{ll}
d\left(\begin{array}{c}
\theta^1 \\
\theta^2 \\
\theta^3\\
\theta^4
\end{array}\right)=\left(\begin{array}{cccc}
\pi^{\prime1} & 0 & 0&0 \\
\pi^{\prime2}&  \pi^{\prime3} & 0&0 \\
\pi^{\prime4}& \pi^{\prime5} & \pi^{\prime6}&0\\
\pi^{\prime7}&\pi^{\prime8}&0&\pi^{\prime9}
\end{array}\right) \wedge\left(\begin{array}{c}
\theta^1 \\
\theta^2 \\
\theta^3\\
\theta^4
\end{array}\right)+\left(\begin{array}{c}
T_{24}^1~ \theta^2 \wedge \theta^4 \\
T_{34}^2~\theta^3\wedge\theta^4 \\
0\\
0
\end{array}\right),
\end{array}
\end{equation}
The essential torsion coefficients are $T_{24}^1=-\frac{a_1}{a_3 a_9}$ and $T_{34}^2=-\frac{a_3}{a_6 a_9}$
which can be normalized  to $-1$ by setting $a_6=\frac{a^2_3} {a_1},a_9=\frac{a_1}{a_3}$. Thus, the structure group is reduced to
\begin{equation}\label{2.8}
{G}_1=\left\{
\left(\begin{array}{cccc}
a_1 & 0 & 0&0 \\
a_2 & a_3 & 0 &0\\
a_4 & a_5 & \frac{a^2_3} {a_1}&0\\
a_7&a_8&0&\frac{a_1}{a_3}
\end{array}\right)
 \Bigg \vert a_1a_3\neq 0
\right\}.
\end{equation}
In the \emph{second loop} through the reduction procedure, the structure equations become, after absorption,
\begin{equation}\label{2.10}
d\left(\begin{array}{c}
\theta^1 \\
\theta^2 \\
\theta^3\\
\theta^4
\end{array}\right)=\left(\begin{array}{cccc}
 \pi^{\prime1} & 0 & 0&0 \\
\pi^{\prime2} &  \pi^{\prime3} & 0&0 \\
\pi^{\prime4} &\pi^{\prime5}& 2\pi^{\prime3}-\pi^{\prime1}&0\\
\pi^{\prime6}&\pi^{\prime7}&0&\pi^{\prime1}-\pi^{\prime3}
\end{array}\right) \wedge\left(\begin{array}{c}
\theta^1 \\
\theta^2 \\
\theta^3\\
\theta^4
\end{array}\right)+\left(\begin{array}{c}
- \theta^2 \wedge \theta^4 \\
 - \theta^3 \wedge \theta^4 \\
T_{34}^3~\theta^3 \wedge \theta^4 \\
0
\end{array}\right).
\end{equation}
The essential torsion coefficient $T_{34}^3=\frac{a_3^2I_1-3a_1a_5+3a_2a_3}{a_1a_3}$
can be normalized to zero by setting\\
  $a_5=\frac{a_2a_3}{a_1}+\frac{a_3^2}{3a_1}I_1$, where $I_1=-f_q$. 
Thus, the structure group is reduced to
\begin{equation}\label{2.8}
{G}_2=\left\{
\left(\begin{array}{cccc}
a_1 & 0 & 0&0 \\
a_2 & a_3 & 0 &0\\
a_4 &\frac{a_2a_3}{a_1}+\frac{a_3^2 }{3a_1}I_1& \frac{a^2_3} {a_1}&0\\
a_7&a_8&0&\frac{a_1}{a_3}
\end{array}\right)
 \Bigg \vert a_1a_3\neq 0
\right\}.
\end{equation}
In the \emph{third loop} through the reduction procedure, the structure equations become, after absorption,
\begin{equation}\label{2.7}
\begin{array}{ll}
d\left(\begin{array}{c}
\theta^1 \\
\theta^2 \\
\theta^3\\
\theta^4
\end{array}\right)=\left(\begin{array}{cccc}
\pi^{\prime1} & 0 & 0&0 \\
\pi^{\prime2}&  \pi^{\prime3} & 0&0 \\
 \pi^{\prime4}& \pi^{\prime2} & 2\pi^{\prime3}-\pi^{\prime1}&0\\
 \pi^{\prime5}&\pi^{\prime6}&0&\pi^{\prime1}-\pi^{\prime3}
\end{array}\right) \wedge\left(\begin{array}{c}
\theta^1 \\
\theta^2 \\
\theta^3\\
\theta^4
\end{array}\right)+\left(\begin{array}{c}
-\theta^2 \wedge \theta^4 \\
-\theta^3\wedge\theta^4 \\
T_{24}^3~\theta^2 \wedge \theta^4\\
0
\end{array}\right).
\end{array}
\end{equation}
The essential torsion coefficient $T_{24}^3=\frac{a_3^2I_2-2a_1a_4+a_2^2}{a_1^2}$
can be normalized to zero by setting\\  $a_4=\frac{a_2^2}{2a_1}+\frac{a_3^2}{2a_1}I_2$, where  $I_2=-\frac{2}{9}I_1^2-f_p-\frac{1}{3}\hat{D}_xI_1$. 
Consequently, the structure group reduces to
\begin{equation}\label{2.8}
{G}_3=\left\{
\left(\begin{array}{cccc}
a_1 & 0 & 0&0 \\
a_2 & a_3 & 0 &0\\
\frac{a_2^2}{2a_1}+\frac{a_3^2}{2a_1}I_2 &\frac{a_2a_3}{a_1}+\frac{a_3^2}{3a_1}I_1 & \frac{a^2_3} {a_1}&0\\
a_7&a_8&0&\frac{a_1}{a_3}
\end{array}\right)
 \Bigg \vert a_1a_3\neq 0
\right\}.
\end{equation}
This yields the adapted coframe (\ref{2.6}) with $S\in G_3$. \\
Proceeding to the \emph{fourth loop} of the reduction algorithm, one finds that, after absorption, the structure equations become
\begin{equation}\label{2.7}
\begin{array}{ll}
d\left(\begin{array}{c}
\theta^1 \\
\theta^2 \\
\theta^3\\
\theta^4
\end{array}\right)=\left(\begin{array}{cccc}
\pi^{\prime1} & 0 & 0&0 \\
\pi^{\prime2}&  \pi^{\prime3} & 0&0 \\
0& \pi^{\prime2} & 2\pi^{\prime3}-\pi^{\prime1}&0\\
\pi^{\prime4}&\pi^{\prime5}&0&\pi^{\prime1}-\pi^{\prime3}
\end{array}\right) \wedge\left(\begin{array}{c}
\theta^1 \\
\theta^2 \\
\theta^3\\
\theta^4
\end{array}\right)+\left(\begin{array}{c}
-\theta^2 \wedge \theta^4 \\
-\theta^3\wedge\theta^4 \\
T_{14}^3~\theta^1 \wedge \theta^4\\
0
\end{array}\right).
\end{array}
\end{equation}
The essential  torsion coefficient is
\begin{equation}
T^3_{14}=\frac{a_3^3}{a_1^3}I_3,\\
\end{equation}
where 
\begin{equation}
I_3=-\frac{1}{3}I_1I_2-f_u-\frac{1}{2}\hat{D}_xI_2. 
\end{equation}
As a result, we obtain the following branch.
\subsection{Branch  $I_3=0$}
In this branch, the structure equations have the form
\begin{equation}\label{3.33}
\begin{array}{ll}
d\left(\begin{array}{c}
\theta^1 \\
\theta^2 \\
\theta^3\\
\theta^4
\end{array}\right)=\left(\begin{array}{cccc}
\pi^{\prime1} & 0 & 0&0 \\
\pi^{\prime2}& \pi^{\prime3}& 0&0 \\
0& \pi^{\prime2}&2\pi^{\prime3}-\pi^{\prime1}&0\\
\pi^{\prime4}&\pi^{\prime5}&0&\pi^{\prime1}-\pi^{\prime3}
\end{array}\right) \wedge\left(\begin{array}{c}
\theta^1 \\
\theta^2 \\
\theta^3\\
\theta^4
\end{array}\right)+
\left(\begin{array}{c}
-\theta^2 \wedge \theta^4 \\
-\theta^3\wedge\theta^4   \\
0 \\
0
\end{array}\right),
\end{array}
\end{equation}
in which there are no group-dependent invariants among the essential torsion coefficients and their coframe derivatives, so Cartan’s test must be applied.

In this loop,   the  Maurer-Cartan forms  $\pi^{\kappa},~\kappa=1, \ldots, 5$,  are given explicitly as
\begin{equation}\label{2.11}
\left(\begin{array}{c}
\pi^{1} \\
\pi^{2} \\
\pi^{3}\\
\pi^{4}\\
\pi^{5}\\
\end{array}\right)=\Pi_{{G_3}}
\left(\begin{array}{c}
da_1 \\
da_2 \\
da_3\\
da_7\\
da_8
\end{array}\right)
+\Pi_M
\left(\begin{array}{c}
dx \\
du \\
dp\\
dq
\end{array}\right),
\end{equation}
where the matrices $\Pi_{{G_3}}$ and $\Pi_M$  are the projections of the forms $\pi^{\kappa},~\kappa=1, \ldots, 5,$ on $T^*{{G_3}}$ and $T^*M$ respectively, and are given by
\begin{equation}\label{2.12}
\Pi_{{G_3}}=
\left(\begin {array}{ccccc}
 \frac{1}{a_1}&0&0&0&0\\ 
0& \frac{1}{a_1}&-\frac{a_2}{a_1a_3}&0&0\\ 
0&0& \frac{1}{a_3}&0&0\\ 
\frac{\Delta}{a_1^2a_3}&0& -\frac{\Delta}{a_1a_3^2}&\frac{1}{a_1}&-\frac{a_2}{a_1a_3}\\
-\frac{a_8}{a_1a_3}&0&\frac{a_8}{a_3^2}&0&\frac{1}{a_3}
\end {array} \right),
\Pi_M=
\left(\begin {array}{cccc}
0&0&0&0\\ 
0&0&0&0\\ 
0&0&0&0\\ 
0&0&0&0\\ 
0&0&0&0\\ 
\end {array} \right),
\end{equation}
where $\Delta=a_2a_8-a_3a_7$.
The modified Maurer-Cartan forms are given by the general solution to the linear absorption system as
\begin{equation}\label{2.13}
\left(\begin{array}{c}
\pi^{\prime1} \\
\pi^{\prime2} \\
\pi^{\prime3}\\
\pi^{\prime4}\\
\pi^{\prime5}\\
\end{array}\right)=
\left(\begin{array}{c}
\pi^{1} \\
\pi^{2} \\
\pi^{3}\\
\pi^{4}\\
\pi^{5}\\
\end{array}\right)
+\left(\Lambda^{(1)}+\Lambda^{(2)}\right)
\left(\begin{array}{c}
\theta^1 \\
\theta^2 \\
\theta^3\\
\theta^4
\end{array}\right),
\end{equation}
where $\Lambda^{(1)}$ is the matrix of determined parameters
\begin{equation}
\Lambda^{(1)}=
\begin{pmatrix}
-\dfrac{A+2B}{6a_1^2a_3}
& \dfrac{a_7}{a_1}
& 0
& -\dfrac{a_2}{a_1}
\\[2.2ex]
\dfrac{H}{6a_1^3a_3}
& -\dfrac{A+B+C}{6a_1^2a_3}
& -\dfrac{\Delta}{a_1a_3}
& \dfrac{Q}{6a_1^2}
\\[2.2ex]
-\dfrac{A+B}{6a_1^2a_3}
&-\dfrac{\Delta}{a_1a_3}+\dfrac{I_1a_8}{3a_1}-\dfrac{I_{1_q}}{6a_3}
& \dfrac{a_8}{a_3}
& -\dfrac{a_3I_1}{3a_1}
\\[2.2ex]
0
& \dfrac{a_8C+6a_3a_7\Delta}{6a_1^2a_3^2}
& -\dfrac{a_8\Delta}{a_1a_3^2}
& 0
\\[2.2ex]
0
& 0
& \dfrac{a_8^2}{a_3^2}
& \frac{a_7}{a_1}
+\frac{I_{1_q}}{6a_3}
-\frac{2I_1a_8}{3a_1}
\end{pmatrix}.
\end{equation}
which is expressed in terms of
\begin{equation*}
\begin{array}{lll}
Q&=2I_1a_2a_3-3I_2a_3^2-3a_2^2,\\
A&=2(2I_1a_3-3a_2)\Delta
   +a_1\bigl(3a_3I_{2_q}-2a_2I_{1_q}\bigr),\\
B&=a_8Q+6a_2\Delta,\quad
C=a_8Q-2I_1a_3\Delta,\\
H&=3a_3^3a_8\hat{D}_x(I_2)+2I_1I_2a_3^3a_8+6a_3^3a_8f_u 
+a_1a_3^2\bigl(I_1I_{2_q}-I_2I_{1_q}-3I_{2_p}+2I_{1_u}\bigr)\\
 &+3a_1a_2a_3I_{2_q}-a_1a_2^2I_{1_q}+\Delta Q,
\end{array}
\end{equation*}
and the matrix  $\Lambda^{(2)}$ contains the indeterminate parameters and is given by
\begin{equation}\label{2.15}
\Lambda^{(2)}=
\begin{pmatrix}
  2r_2&0&0&0 \\
  0&r_2&0&0\\
  r_2&0&0&0 \\
 r_1&r_3&0&r_2 \\
  r_3&r_4&0&0 \\
\end{pmatrix}.
\end{equation}

The general solution to the linear absorption system has four free parameters, and  the degree of indeterminacy is $r^{(1)}=4$. 

Moreover, the coefficient matrix $L^i_\rho[v]$ for the Maurer–Cartan forms $\pi^{'\rho}$ in the structure equation (\ref{3.33}) upon replacing each coframe element $\theta^i$ by the corresponding entry $v_i$ of the vector $v$ (and deleting the wedge products) can be given as
\begin{equation}\label{2.16}
L[v]=\begin{pmatrix}
v^{1}&0&0&0&0\\
0&v^{1}&v^{2}&0&0\\
-v^{3}&v^{2}&2v^{3}&0&0\\
v^{4}&0&-v^{4}&v^{1}&v^{2}
\end{pmatrix}
\end{equation}
Accordingly, the Cartan characters $s^{\prime}_{1}, s^{\prime}_{2}, s^{\prime}_{3}, s^{\prime}_{4}$ for the lifted coframe
can be determined as follows:
\begin{equation*}
\begin{aligned}
s^{\prime}_1&=max\{rank ~ L[v] \mid v\in \mathbb{R}^4\},\\
s^{\prime}_1+s^{\prime}_2&=max\bigg\{rank\begin{pmatrix}L[v_1]\\
L[v_2]\end{pmatrix}\mid v_1,v_2\in\mathbb{R}^4\bigg\},\\
s^{\prime}_1+s^{\prime}_2+s^{\prime}_3&=max\Biggl\{rank\begin{pmatrix}L[v_1]\\
L[v_2]\\
L[v_3]\\
\end{pmatrix}\mid v_1, v_2, v_3\in\mathbb{R}^4\Biggl\},\\
s^{\prime}_1+s^{\prime}_2+s^{\prime}_3+s^{\prime}_4&=r,
\end{aligned}
\end{equation*}
where $r$ is the dimension of the structure group $G_3$.
Hence, the Cartan characters are $s^{\prime}_{1}=4, s^{\prime}_{2}=1, s^{\prime}_{3}=0$ and $s^{\prime}_{4}=0$. Therefore, $s^{\prime}_{1}+2s^{\prime}_{2}+3s^{\prime}_{3}+4s^{\prime}_{4}= 6 > r^{(1)} = 4$, and the Cartan involutivity test does not
hold. Hence, the lifted coframe fails Cartan's involutivity criterion and prolongation is  required.

The prolonged structure group is a $\textsl{four-dimensional
abelian group}$ having the $9\times9$ matrix representation
\begin{equation}\label{b21}
G^{(1)}=\left\{ \left(
\begin{array}{cc}
  I & 0 \\
  \Lambda^{(2)} & I \\
\end{array}%
\right) \,\middle\vert\,
\Lambda^{(2)}=\left(%
\begin{array}{cccc}
  2r_2&0&0&0 \\
  0&r_2&0&0\\
  r_2&0&0&0 \\
 r_1&r_3&0&r_2 \\
  r_3&r_4&0&0 \\
\end{array}%
\right) \right \},
\end{equation}
and the prolonged base coframe on the space $M^{(1)}=M\times {G_3}$ consists of the forms (\ref{2.6}) with $S\in G_3$ and the modified Maurer-Cartan forms (\ref{2.13}) after setting $\Lambda^{(2)}$ equal to the zero matrix.

The resulting structure equations associated with the lifted prolonged coframe in the absorbed form are the original  structure equations (\ref{3.33}) and
\begin{equation}\label{3.34}
\begin{array}{lll}
&d\pi^{\prime1}=2\rho^{\prime2}\wedge\theta^{1}-\theta^{2} \wedge \pi^{\prime4}+\theta^{4}\wedge \pi^{\prime2},\\
&d\pi^{\prime2}=\rho^{\prime2}\wedge\theta^{2}+T^6_{13}~\theta^{1}\wedge\theta^{3}+T^6_{23}~\theta^{2}\wedge\theta^{3}-\theta^{3}\wedge\pi^{\prime4}-\pi^{\prime1}\wedge\pi^{\prime2}-\pi^{\prime2}\wedge\pi^{\prime3},\\
&d\pi^{\prime3}=\rho^{\prime2}\wedge\theta^{1}+T^7_{12}~\theta^1\wedge\theta^2+T^7_{13}~\theta^1\wedge\theta^3+T^7_{23}~\theta^2\wedge\theta^3-\theta^2\wedge\pi^{\prime4}-\theta^3\wedge\pi^{\prime5},\\
&d\pi^{\prime4}=\rho^{\prime1}\wedge\theta^1+\rho^{\prime3}\wedge\theta^2+\rho^{\prime2}\wedge\theta^4+T^8_{34}~\theta^{3}\wedge\theta^{4}-\pi^{\prime2}\wedge\pi^{\prime5}-\pi^{\prime3}\wedge\pi^{\prime4},\\
&d\pi^{\prime5}=\rho^{\prime3}\wedge\theta^1+\rho^{\prime4}\wedge\theta^2+T^9_{14}~\theta^{1}\wedge\theta^{4}+T^9_{34}~\theta^{3}\wedge\theta^{4}-\theta^4\wedge\pi^{\prime4}+\pi^{\prime1}\wedge\pi^{\prime5}-2\pi^{\prime3}\wedge\pi^{\prime5},
\end{array}
\end{equation}
where $\rho^{\prime\kappa},~\kappa=1, \ldots, 4,$ are the modified Maurer-Cartan forms.\\
We normalize  $T^{6}_{13}=0, T^{6}_{23}=0, T^7_{23}=0$ by setting
\begin{equation}\label{2.25}
\begin{array}{ll}
r_1 &= \,\frac{{a_2 ^2 a_7 a_8 }}{{2\,a_1 ^3 a_3 }} - \frac{{a_2 ^3 a_8 ^2 }}{{2\,a_1 ^3 a_3 ^2 }} - \left( {\frac{{a_2 a_7 a_8 }}{{a_1 ^3 }} - \,\frac{{2\,a_2 ^2 a_8 ^2 }}{{3\,a_1 ^3 a_3 }} - \,\frac{{a_3 a_7 ^2 }}{{3\,a_1 ^3 }}} \right)I_1  - \,\left( {\frac{{a_2 a_8 ^2 }}{{2\,a_1 ^3 }} - \,\frac{{a_3 a_7 a_8 }}{{2\,a_1 ^3 }}} \right)I_2 \\
 &- \, \left( {\frac{{a_2 ^2 a_8 }}{{2\,a_1 ^2 a_3 ^2 }} - \, \frac{{a_2 a_7 }}{{3\,a_1 ^2 a_3 }}} \right) I_{1_q} - \left( {\frac{{a_7 }}{{2\,a_1 ^2 }} - \,\frac{{a_2 a_8 }}{{\,a_1 ^2 a_3 }}} \right) I_{2_q}+ \, \frac{{a_2 ^2 }}{{6\,a_1 a_3 ^3 }} I_{1_{qq}}- \,\frac{{a_2 }}{{2\,a_1 a_3 ^2 }} I_{2_{qq}}\\
 &+ \,\frac{{a_8 }}{{6\,a_1 ^2 }}\left( {-3 I_{2_p} + 2 I_{1_u}+  I_1 I_{2_q} - I_2 I_{1_q}} \right) + \,\frac{1}{{6\,a_1 a_3 }}\left( { I_2 I_{1_{qq}} -  I_1 I_{2_{qq}} - 2 I_{1_{uq}} + 3 I_{2_{uq}}} \right), \\
r_3 &= \, - \,\frac{{a_7 ^2 }}{{2\,a_1 ^2 }} + \,\frac{{a_2 ^2 a_8 ^2 }}{{2\,a_1 ^2 a_3 ^2 }} - \,\left( {\frac{{2\,a_2 a_8 ^2 }}{{3\,a_1 ^2 a_3 }} - \,\frac{{2\,a_7 a_8 }}{{3\,a_1 ^2 }}} \right)I_1  - \left( {\frac{{a_7 }}{{6\,a_1 a_3 }} - \frac{{a_2 a_8 }}{{2\,a_1 a_3 ^2 }}} \right) I_{1_q} - \frac{{a_8 }}{{2\,a_1 a_3 }} I_{2_q} \\
&- \frac{{a_2 }}{{6\,a_3 ^3 }}I_{1_{qq}}+ \frac{1}{{4\,a_3 ^2 }} I_{2_{qq}}, \\
r_4 &=  - \frac{{a_2 a_8 ^2 }}{{a_1 a_3 ^2 }} + \,\frac{{2\,a_8 ^2 }}{{3\,a_1 a_3 }}I_1  - \frac{{a_8 }}{{2a_3 ^2 }}I_{1_q}+ \frac{{a_1 }}{{6\,a_3 ^3 }}I_{1_{qq}}. \\
\end{array}
\end{equation}	
These normalizations  reduce the structure group to
 the {\it one-dimensional group} $G^{(1)}_1$
\begin{equation}
G^{(1)}_1=\left\{ \left(
\begin{array}{cc}
  I & 0 \\
R& I \\
\end{array}%
\right) \right\}
\end{equation}
where $R$ is the matrix $\Lambda^{(2)} $ given in (\ref{2.15}) after
incorporating the values of the parameters $r_1, r_3, r_4$ obtained in (\ref{2.25}). This leads to the structure equation, after absorption, given by the original  structure equations (\ref{3.33}) and
\begin{equation}\label{3.35}
\begin{array}{lll}
&d\pi^{\prime1}=2\rho^{\prime2}\wedge\theta^{1}-\theta^2\wedge\pi^{\prime4}+\theta^4\wedge\pi^{\prime2},\\
&d\pi^{\prime2}=\rho^{\prime2}\wedge\theta^{2}-\theta^3\wedge\pi^{\prime4}-\pi^{\prime1}\wedge\pi^{\prime2}-\pi^{\prime2}\wedge\pi^{\prime3},\\
&d\pi^{\prime3}=\rho^{\prime2}\wedge\theta^{1}-\theta^2\wedge\pi^{\prime4}-\theta^3\wedge\pi^{\prime5},\\
&d\pi^{\prime4}=\rho^{\prime2}\wedge\theta^4+T^8_{12}~\theta^1\wedge\theta^2+T^8_{13}~\theta^1\wedge\theta^3+T^8_{23}~\theta^2\wedge\theta^3-\pi^{\prime2}\wedge\pi^{\prime5}-\pi^{\prime3}\wedge\pi^{\prime4},\\
&d\pi^{\prime5}=T^9_{12}~\theta^{1}\wedge\theta^{2}+T^9_{13}~\theta^{1}\wedge\theta^{3}+T^9_{23}~\theta^{2}\wedge\theta^{3}-\theta^4\wedge\pi^{\prime4}+\pi^{\prime1}\wedge\pi^{\prime5}-2\pi^{\prime3}\wedge\pi^{\prime5},
\end{array}
\end{equation}
The essential torsion coefficient $T^9_{23}=-\frac{a_1^2I_4}{6a_3^5}$, where $I_4=I_{1_{qqq}}$ is a relative invariant.
The relative invariant $I_4=0$  for the canonical form
$u''' = 0$. Thus, we choose the following branch.
\subsubsection {Branch  $I_4=0$}
In this branch,  there is no further unabsorbable torsion left. Moreover, $\rho^{\prime2}$
is uniquely determined; hence, the equivalence problem is fully determined. This yields an $e$-structure on the ten-dimensional prolonged space $M^{(2)}=M^{(1)} \times G^{(1)}_1$. The resulting coframe consists of the original lifted coframe \begin{equation}\label{2.27}
\begin{array}{ll}
\left(\begin{array}{c}
\theta^1 \\
\theta^2 \\
\theta^3\\
\theta^4
\end{array}\right)=\left(\begin{array}{cccc}
a_1 & 0 & 0&0 \\
a_2 & a_3 & 0 &0\\
\frac{a_2^2}{2a_1}+\frac{a_3^2}{2a_1}I_2 &\frac{a_2a_3}{a_1}+\frac{a_3^2}{3a_1}I_1 & \frac{a^2_3} {a_1}&0\\
a_7&a_8&0&\frac{a_1}{a_3}
\end{array}\right)\left(\begin{array}{l}
\omega^1\\\omega^2\\\omega^3\\\omega^4\\
\end{array}\right),
\end{array}
\end{equation}
\begin{equation}\label{2.28}
\left(\begin{array}{c}
\pi^{\prime1} \\
\pi^{\prime2} \\
\pi^{\prime3}\\
\pi^{\prime4}\\
\pi^{\prime5}\\
\end{array}\right)=
\left(\begin{array}{c}
\pi^{1} \\
\pi^{2} \\
\pi^{3}\\
\pi^{4}\\
\pi^{5}\\
\end{array}\right)
+\left(\Lambda^{(1)}+R\right)
\left(\begin{array}{c}
\theta^1 \\
\theta^2 \\
\theta^3\\
\theta^4
\end{array}\right),
\end{equation}
together with the modified Maurer-Cartan form $\rho^{\prime 2}$. 
This results in the  structure equations 
\begin{equation}\label{2.277}
\begin{array}{lll}
d\theta^1&=&-\theta^1\wedge\pi^{\prime1}
            -\theta^2\wedge\theta^4,\\[1ex]

d\theta^2&=&-\theta^1\wedge\pi^{\prime2}
            -\theta^2\wedge\pi^{\prime3}
            -\theta^3\wedge\theta^4,\\[1ex]

d\theta^3&=&-\theta^2\wedge\pi^{\prime2}
            +\theta^3\wedge\pi^{\prime1}
            -2\theta^3\wedge\pi^{\prime3},\\[1ex]

d\theta^4&=&-\theta^1\wedge\pi^{\prime4}
            -\theta^2\wedge\pi^{\prime5}
            -\theta^4\wedge\pi^{\prime1}
            +\theta^4\wedge\pi^{\prime3},\\[1ex]

d\pi^{\prime1}&=&2\rho^{\prime2}\wedge\theta^1
                 -\theta^2\wedge\pi^{\prime4}
                 +\theta^4\wedge\pi^{\prime2},\\[1ex]

d\pi^{\prime2}&=&\rho^{\prime2}\wedge\theta^2
                 -\theta^3\wedge\pi^{\prime4}
                 -\pi^{\prime1}\wedge\pi^{\prime2}
                 -\pi^{\prime2}\wedge\pi^{\prime3},\\[1ex]

d\pi^{\prime3}&=&\rho^{\prime2}\wedge\theta^1
                 -\theta^2\wedge\pi^{\prime4}
                 -\theta^3\wedge\pi^{\prime5},\\[1ex]

d\pi^{\prime4}&=&\rho^{\prime2}\wedge\theta^4
                 -\pi^{\prime2}\wedge\pi^{\prime5}
                 -\pi^{\prime3}\wedge\pi^{\prime4},\\[1ex]

d\pi^{\prime5}&=&-\theta^4\wedge\pi^{\prime4}
                 +\pi^{\prime1}\wedge\pi^{\prime5}
                 -2\pi^{\prime3}\wedge\pi^{\prime5},\\[1ex]

d\rho^{\prime2}& =&-\pi^{\prime2}\wedge\pi^{\prime4}
-\pi^{\prime1}\wedge\rho^{\prime2}.
\end{array}
\end{equation}
The invariant structure of the prolonged coframe are all constant.
We have produced an invariant coframe with rank zero on the
ten-dimensional space with coordinates $x,u,p, q, a_1, a_2,$ $a_3, a_7,
a_8, r_2$. Any such differential equation admits a
ten-dimensional symmetry group of contact transformations. This proves the following theorem.
\begin{theorem} \label{Th}
A scalar third-order ODE of the form $u''' = f(x,u,u',u'')$ is equivalent under a contact
transformation (\ref{2.4}) to the canonical form $\bar{u}'''=0$, and hence admits
the maximal ten-dimensional contact symmetry Lie algebra, if and only if
the relative invariants $I_3$ and $I_4$ vanish identically.
\end{theorem}
\begin{remark}
The modified Maurer-Cartan form $\rho^{\prime 2}$ has a complicated expression; however, using (\ref{2.277}), it can be expressed compactly through its contraction with the invariant differential operator $\partial\theta^4$ as follows
\begin{equation}\label{2.28}
\rho'^2 = -\frac{\partial}{\partial \theta^4} \righthalfcup d\pi'^4
\end{equation}
where $\partial\theta^4$ is  given by 
\begin{equation}\label{2.29}
\begin{array}{lll}
\partial\theta^4={}&
\frac{a_3}{a_1}\hat{D}_x
+a_2\frac{\partial}{\partial a_1}
+\left(\frac{a_3^2}{2a_1}I_2+\frac{a_2^2}{2a_1}\right)\frac{\partial}{\partial a_2}
+\frac{a_3^2}{3a_1}I_1\frac{\partial}{\partial a_3}
+\left(
\frac{2a_2a_8}{3a_1}I_1
-\frac{a_3a_7}{3a_1}I_1
-a_1r_2
-\frac{a_2}{6a_3}I_{1_q}\right)
\frac{\partial}{\partial a_7}\\
&+\left(
\frac{a_3a_8}{3a_1}I_1
-\frac{1}{6}I_{1_q}
+\frac{a_2a_8}{a_1}
-\frac{a_3a_7}{a_1}\right)
\frac{\partial}{\partial a_8}
+(\frac{a_3  r_2}{3a_1}I_1
+\frac{a_2^2a_7}{2a_1^3}
-\frac{a_2^2a_8 }{3a_1^3}I_1
+\frac{a_2^2 }{12a_1^2a_3}I_{1_q}
-\frac{a_3^2a_8}{2a_1^3}\hat{D}_x I_2\\
&+
\frac{a_3^2a_7}{a_1^3}
\left(
\frac{1}{2}I_2-\frac{1}{3}\hat{D}_x I_1
\right)
+\frac{2a_2a_3a_8}{3a_1^3}\hat{D}_x I_1
-\frac{a_3}{a_1^2}I_5
+\frac{a_2}{a_1^2}I_6)
\frac{\partial}{\partial r_2},
\end{array}
\end{equation}
which is given in terms of
\begin{equation}
\begin{array}{ll}
I_5  = -\frac{2}{{3\,}}I_{1_u} +\,\frac{1}{2}I_{2_p} - \,\frac{1}{{6\,}}I_1 I_{2_q}+ \,\frac{1}{{4\,}}I_2 I_{1_q},
&I_6  = \,-\frac{1}{{3\,}}I_{1_p} + \frac{1}{{2\,}}I_{2_q} + \,\frac{1}{{18\,}}I_1 I_{1_q}. 
\end{array}
\end{equation}
\end{remark}
\begin{remark}
The six modified Maurer–Cartan forms can be written in matrix form as
\begin{equation}\label{2.30}
\begin{pmatrix}
\pi^{\prime 1}\\
\pi^{\prime 2}\\
\pi^{\prime 3}\\
\pi^{\prime 4}\\
\pi^{\prime 5}\\
\rho^{\prime2}
\end{pmatrix}
=
\left(
\begin{array}{cc}
\Pi_{G_3} &0\\
\boldsymbol{\beta} & 1
\end{array}
\right)
\begin{pmatrix}
da_1\\
da_2\\
da_3\\
da_7\\
da_8\\
dr_2
\end{pmatrix}+
\begin{pmatrix}
\Lambda^{(1)}+R
\\[2mm]
s_1 \quad s_2 \quad s_3 \quad s_4
\end{pmatrix}
\begin{pmatrix}
\theta^1\\
\theta^2\\
\theta^3\\
\theta^4
\end{pmatrix},
\end{equation}
where $\boldsymbol{\beta}^{T}
=
\begin{pmatrix}
\displaystyle
\frac{r_2}{a_1}
-\frac{a_2^2a_8}{2a_1^3a_3}
+\left(
-\frac{a_3a_7}{3a_1^3}
+\frac{2a_2a_8}{3a_1^3}
\right)I_1
-\frac{a_3a_8}{2a_1^3}I_2\\
\displaystyle
\frac{a_7}{a_1^2}
-\frac{2a_8}{3a_1^2}I_1
+\frac{I_{1_q}}{6a_1a_3}\\
\displaystyle
\frac{a_2^2a_8}{2a_1^2a_3^2}
-\frac{a_2a_7}{a_1^2a_3}
+\frac{a_7}{3a_1^2}I_1
+\frac{a_8}{2a_1^2}I_2
-\frac{a_2}{6a_1a_3^2}I_{1_q}\\
\displaystyle
\frac{a_3}{3a_1^2}I_1\\
\displaystyle
\frac{a_2^2}{2a_1^2a_3}
-\frac{2a_2}{3a_1^2}I_1
+\frac{a_3}{2a_1^2}I_2
\end{pmatrix}$, \\
and $s_i$ denotes the coefficient of $\rho^{\prime 2}$ in the $\theta^i$-direction, for $i=1,\ldots,4$.\\
Substituting the contact transformation (\ref{2.4}) into the symmetrized Cartan formulation of (\ref{2.30}), together with
$\bar{f}=0$, for which $\bar{I}_1=\bar{I}_2=0$, and
$\bar{a}_1=1$, $\bar{a}_2=0$, $\bar{a}_3=1$,
$\bar{a}_7=0$, $\bar{a}_8=0$, $\bar{r}_2=0$, gives
$\bar{s}_i=0$, $i=1,\ldots,4$, 
$\bar{\rho}^{\prime2}=0$ and
$\bar{\pi}^{\prime i}=0$, $i=1,\ldots,5$, yield the following first-order system of PDEs:
\begin{equation}\label{2.31}
\begin{pmatrix}
\mathcal D_1(a_1) & \mathcal D_2(a_1) &
\mathcal D_3(a_1) & \mathcal D_4(a_1)
\\
\mathcal D_1(a_2) & \mathcal D_2(a_2) &
\mathcal D_3(a_2) & \mathcal D_4(a_2)
\\
\mathcal D_1(a_3) & \mathcal D_2(a_3) &
\mathcal D_3(a_3) & \mathcal D_4(a_3)
\\
\mathcal D_1(a_7) & \mathcal D_2(a_7) &
\mathcal D_3(a_7) & \mathcal D_4(a_7)
\\
\mathcal D_1(a_8) & \mathcal D_2(a_8) &
\mathcal D_3(a_8) & \mathcal D_4(a_8)
\\
\mathcal D_1(r_2) & \mathcal D_2(r_2) &
\mathcal D_3(r_2) & \mathcal D_4(r_2)
\end{pmatrix}=
-
\left(
\begin{array}{cc}
\Pi_{G_3} & 0\\
\boldsymbol{\beta} & 1
\end{array}
\right)^{-1}
\begin{pmatrix}
\Lambda^{(1)}+R
\\[1mm]
s_1\;\;s_2\;\;s_3\;\;s_4
\end{pmatrix},
\end{equation}
where the  differential operators 
\(\mathcal D_1,\mathcal D_2,\mathcal D_3,\mathcal D_4\) are dual to the lifted coframe (\ref{2.6}) for $S\in G_3$ on $M=\mathbf{J}^2$, and can be obtained as follows:
\begin{align*}
\mathcal D_1
&=
\frac{\Delta}{a_1^2}\hat{D}_x
+\frac{1}{a_1}\partial_u
-\frac{a_2}{a_1a_3}\partial_p
+\frac{2a_2a_3I_1-3a_3^2I_2+3a_2^2}
       {6a_1a_3^2}\partial_q,\\[1mm]
\mathcal D_2
&=
-\frac{a_8}{a_1}\hat{D}_x
+\frac{1}{a_3}\partial_p
-\frac{a_3I_1+3a_2}{3a_3^2}\partial_q,\\[1mm]
\mathcal D_3
&=\frac{a_1}{a_3^2}\partial_q,
\qquad
\mathcal D_4=\frac{a_3}{a_1}\hat{D}_x.
\end{align*}
\end{remark}
\section{Construction of Contact Transformation Using the Invariant Coframe}
\setcounter{equation}{0}
We now use the invariant coframe to construct a contact transformation between $u'''=f(x,u,p,q)$ and $\bar{u}'''=0$. The construction is summarized in the following theorem.
\begin{theorem} 
Assume that the third-order ODEs
\begin{equation}\label{eee4.8}
\begin{array}{lll}
u^{\prime \prime\prime}=f\left(x, u, u^{\prime},u^{\prime \prime}\right), \quad \bar{u}^{\prime \prime\prime}=0,
\end{array}
\end{equation}
are  equivalent  under the contact transformation
\begin{equation}\label{eee4.9}
\bar x = \varphi (x,u,p),\,\,\bar u = \psi (x,u,p),\,\, \bar{p}=\chi(x,u,p),
\end{equation} 
with a nonvanishing Jacobian. Then the linearizing contact transformation (\ref{eee4.9}) between  the third-order ODEs (\ref{eee4.8}) can be obtained in the following systematic way:
\begin{description} 
\item[Step 1.] Find a solution $(v,w)$, with $v\neq 0$, of the
\textbf{first-order linear system of PDEs}
\begin{equation}\label{4.21}
\widehat{D}_xv=w,\qquad
\widehat{D}_xw+\frac{1}{3}I_1w=-\frac{1}{4}I_2v.
\end{equation} 
\item[Step 2.] Find a nonzero solution $a_3$ of the \textbf{first-order linear system of PDEs}
\begin{equation}\label{3.4}
\begin{array}{l}
\hat{D}_xa_3 =\frac{1}{3}I_1 a_3,\\
v a_{3_p}+(2w-\frac{1}{3}v I_1) a_{3_q} =(\frac{1}{6}vI_{1_q}-2v_p)a_3.
\end{array}
\end{equation}
\item[Step 3.] Set $a_1=\frac{1}{v^2}$, $a_2=-2a_3\frac{w}{v}$.
\item[Step 4.] Set $a_8=\left(\frac{a_1}{a_3}\right)_q$, $a_7=\left(\frac{a_1}{a_3}\right)_p-\hat{D}_xa_8.$
\item[Step 5.] Set $r_2=(\frac{a_3}{a_1^2})(-\frac{1}{3}I_1 a_7+\frac{2}{3}\frac{a_2a_8}{a_3}I_1-\frac{1}{6}\frac{a_1a_2}{a_3^2} I_{1_q}-\hat{D}_xa_7 )$.
\item[Step 6.] Verify, by direct substitution, that the functions obtained in Steps 1–5 satisfy  

\begin{equation}\label{4.21*}
\begin{array}{l}
\hat{D}_xr_2=\frac{1}{3}I_1~r_2+ \,\frac{{a_2 ^2 a_7 }}{{2a_1 ^2 a_3 }} - \frac{{a_2 ^2 a_8 }}{{3\,a_1 ^2 a_3 }}I_1  + \frac{{a_2 ^2 }}{{12\,a_1 a_3 ^2 }}I_{1_q}- \frac{1}{2}\frac{{a_3 a_8 }}{{a_1 ^2 }}\hat{D}_x I_2 \\
 + \frac{{a_3 a_7 }}{{a_1 ^2 }}\left(\frac{1}{2} I_2-\frac{1}{3}\hat{D}_x I_1\right) + \frac{2}{3}\frac{{\,a_2 a_8 }}{{a_1 ^2 }}\hat{D}_x I_1 - \frac{1}{{\,a_1 }}I_5  + \frac{{a_2 }}{{a_1 a_3}}I_6,\\
\end{array}
\end{equation}
and all remaining equations of the full auxiliary system (\ref{2.31}). 
\item[Step 7.] Find a solution $\eta$ of the \textbf{first-order linear system of PDEs}
\begin{equation}\label{4.22}
\begin{array}{l}
\hat{D}_x\eta=0 ,\,\, \eta_u=\frac{a_2^2}{2a_1}+\frac{a_3^2}{2a_1}I_2, \,\,\eta_p=\frac{a_2a_3}{a_1}+\frac{a_3^2}{2a_1}I_1,\,\,\eta_q=\frac{a_3^2}{a_1}.\\
\end{array}
\end{equation}
\item[Step 8.] Find a solution $\chi$ of the \textbf{first-order linear system of PDEs}
\begin{equation}\label{4.23}
\begin{array}{l}
\hat{D}_x\chi=\frac{a_1}{a_3}~\eta, \,\, \chi_u=a_2+a_7~\eta,\,\,\chi_p=a_3+a_8~\eta.\\
\end{array}
\end{equation}
\item[Step 9.] Find a solution $\varphi$ of the \textbf{first-order linear system of PDEs}
\begin{equation}\label{4.24}
\begin{array}{l}
\hat{D}_x \varphi=\frac{a_1}{a_3} ,\,\,\varphi _u=a_7 ,\,\,\varphi_p=a_8.\\
\end{array}
\end{equation}
\item[Step 10.] Find a solution $\psi$ of the \textbf{first-order linear system of PDEs}
\begin{equation}\label{4.25}
\begin{array}{l}
{a_7}\hat{D}_x\psi=\frac{a_1}{a_3}(\psi_u-a_1),\,\,\psi_x=\chi\varphi_x-pa_1,\,\,\psi_p=a_8~\chi.
\end{array}
\end{equation}
\end{description} 
Where $a_1(x,u,p),~a_2(x,u,p,q),~a_3(x,u,p,q),~a_7(x,u,p),~a_8(x,u,p),~v(x,u,p)$,\\
$~w(x,u,p,q)$, and $r_2(x,u,p,q)$ are auxiliary functions, $\chi=\frac{\hat{D}_x\psi}{\hat{D}_x\varphi },\,\, \eta=\frac{\hat{D}_x\chi}{\hat{D}_x\varphi }$, and 
$
I_1=-f_q,
I_2=-\frac{2}{9}I_1^2-f_p-\frac{1}{3}\hat{D}_xI_1.
$
\proof
Inserting the contact transformation (\ref{eee4.9}) into the
symmetrized Cartan formulation of (\ref{2.27}),  and using
$\bar{I}_1=\bar{I}_2=0$ for $\bar{f}=0$ and
$\bar{a}_1=1,\bar{a}_2=0,\bar{a}_3=1,\bar{a}_7=0$, $\bar{a}_8=0$,  results in the following equations 
\begin{equation}
\begin{array}{lll}
\hat{D}_x\eta=0 ,\,\, \eta_u=\frac{a_2^2}{2a_1}+\frac{a_3^2}{2a_1}I_2, \,\,\eta_p=\frac{a_2a_3}{a_1}+\frac{a_3^2}{2a_1}I_1,\,\,\eta_q=\frac{a_3^2}{a_1},\\
\hat{D}_x\chi=\frac{a_1}{a_3}~\eta, \,\, \chi_u=a_2+a_7~\eta,\,\,\chi_p=a_3+a_8~\eta,\\
\hat{D}_x \varphi=\frac{a_1}{a_3} ,\,\,\varphi _u=a_7 ,\,\,\varphi_p=a_8,\\
{a_7}\hat{D}_x\psi=\frac{a_1}{a_3}(\psi_u-a_1),\,\,\psi_x=\chi\varphi_x-pa_1,\,\,\psi_p=a_8~\chi.\\
\end{array}
\end{equation}
Similarly, inserting the contact transformation (\ref{eee4.9}) into the
symmetrized Cartan formulation of (\ref{2.30}), yields the first-order system of PDEs (\ref{2.31}). Consequently, the auxiliary functions are determined by the following linear and Riccati-type differential equations along the contact direction $\theta^4$.
\begin{equation}\label{3.12}
\begin{array}{lll}
\hat{D}_xa_3=\frac{1}{3}I_1a_3,\\
\hat{D}_xa_2=\frac{1}{2a_3}a_2^2+\frac{1}{2}a_3I_2,\\
\hat{D}_xa_1=\frac{a_2}{a_3}a_1,\\
\hat{D}_xa_7 =-\frac{1}{3}I_1 a_7+\frac{2}{3}\frac{a_2a_8}{a_3}I_1-\frac{1}{6}\frac{a_1a_2}{a_3^2} I_{1_q}-\frac{a_1^2}{a_3}r_2,\\
\hat{D}_xa_8 =(\frac{a_2}{a_3}+\frac{1}{3}I_1)a_8-\frac{1}{6}\frac{a_1}{a_3} I_{1_q}-a_7.\\
\hat{D}_xr_2=\frac{1}{3}I_1~r_2+ \,\frac{{a_2 ^2 a_7 }}{{2a_1 ^2 a_3 }} - \frac{{a_2 ^2 a_8 }}{{3\,a_1 ^2 a_3 }}I_1  + \frac{{a_2 ^2 }}{{12\,a_1 a_3 ^2 }}I_{1_q} - \frac{1}{2}\frac{{a_3 a_8 }}{{a_1 ^2 }}\hat{D}_x I_2 \\
 + \frac{{a_3 a_7 }}{{a_1 ^2 }}\left(\frac{1}{2} I_2-\frac{1}{3}\hat{D}_x I_1\right) + \frac{2}{3}\frac{{\,a_2 a_8 }}{{a_1 ^2 }}\hat{D}_x I_1 - \frac{1}{{\,a_1 }}I_5  + \frac{{a_2 }}{{a_1 a_3}}I_6.
\end{array}
\end{equation}
Now, let $a_2=-2a_3\dfrac{w}{v}$, where $w=\hat{D}_xv$. Then the equation
$\hat{D}_x a_2=\dfrac{1}{2a_3}a_2^2+\dfrac{1}{2}a_3I_2$
leads to the system $(\ref{4.21})$. 

Moreover, the equation
$\hat{D}_x a_1=\dfrac{a_2}{a_3}a_1$
gives
$a_1=\dfrac{1}{v^2}$.

Finally, the equations
$\hat{D}_x\varphi=\dfrac{a_1}{a_3}$ and
$\hat{D}_xa_8=\left(\dfrac{a_2}{a_3}+\dfrac{1}{3}I_1\right)a_8-\dfrac{1}{6}\dfrac{a_1}{a_3}I_{1_q}-a_7$
yield\\
$a_8=\left(\dfrac{a_1}{a_3}\right)_q$,
$a_7=\left(\dfrac{a_1}{a_3}\right)_p-\hat{D}_xa_8$, and $v a_{3_p}+\left(2w-\frac{1}{3}vI_1\right)a_{3_q}=\left(\frac{1}{6}vI_{1_q}-2v_p\right)a_3.$
\end{theorem} 
Therefore, Theorem 3.1 provides a constructive algorithm for obtaining the linearizing contact transformation, subject to the verification condition in Step 6, whenever the relative invariants satisfy the conditions of Theorem 2.1.
\begin{remark}
Since $w=\hat{D}_xv$, the function $w$ depends linearly on $q$. Therefore, it can be expressed as
$
w=B(x,u,p)\,q+K(x,u,p),
$
 for some smooth functions $B(x,u,p)$ and $K(x,u,p)$.
\end{remark}
\begin{remark}
The invariant conditions stated in Theorem 2.1 coincide with those previously obtained by Chern \cite{che} and by Neut and Petitot \cite{neu}. The novelty of the present work lies not in the invariant characterization itself, but rather in the explicit Cartan framework developed to construct the corresponding contact transformation through auxiliary functions satisfying only linear and Riccati-type differential equations along the contact direction $\theta^4$.
\end{remark}
\section{Illustrations of the Main Results}
\setcounter{equation}{0}
In this section, we present examples to illustrate our main results.

\begin{example}[A Nonlinear Jerk Equation with Cubic Acceleration Feedback]
In the context of mechanics, \(u(t)\) denotes the position of a particle, \(u''(t)\) is its acceleration, and \(u'''(t)\) is the jerk. The equation
\begin{equation}\label{n2}
u'''=u''^3\\
\end{equation}
may be interpreted as a jerk model in which the rate of change of acceleration is proportional to the cube of the acceleration. As the magnitude of the acceleration increases, the jerk grows much more rapidly, producing abrupt changes in the acceleration and consequently non-uniform motion.

One can verify that the function
$f(x,u,p,q)=q^3$ satisfies the conditions of  Theorem \ref{Th}.
Consequently, this equation admits a ten-dimensional contact symmetry Lie algebra. Moreover, it is equivalent to the canonical form
$\bar{u}'''=0$.
Using $I_1=-3q^2,~I_2=0$ for $u'''=u''^3$, we outline the steps.
\begin{description}
\item[Step 1.]  A solution of system (\ref{4.21}) is $v=1, w=0$. 
 \item[Step 2.]  A solution of system (\ref{3.4}) is $a_3=\frac{1}{q}.$
\item[Step 3.]  The auxiliary functions are $a_1=1$ and $a_2=0.$
\item[Step 4.] The auxiliary functions are $a_7=0$ and $ a_8=1$.
\item[Step 5.] The auxiliary function is $r_2=0$.
\item[Step 6.] The auxiliary functions satisfy the compatibility equation (\ref{4.21*}) and the remaining equations of the
original auxiliary system (\ref{2.31}).
\item[Step 7.] A solution of system (\ref{4.22}) is $\eta=-p-\frac{1}{q}$.
\item[Step 8.] A solution of  system (\ref{4.23}) is $\chi=-x-\frac{1}{2}p^2$.
\item[Step 9.] A solution of  system (\ref{4.24}) is $\varphi =p$.
\item[Step 10.] A solution of  system (\ref{4.25}) is $\psi=-\frac{1}{6}p^3+u-xp$.
\end{description}
Finally, one can verify that the \textbf{proper contact transformation} 
\begin{equation}\label{n4}
\begin{array}{lll}
\varphi \left( x,u,p \right)=p,&\psi \left( x,u,p \right)=-\frac{1}{6}p^3+u-xp,&\chi\left( x,u,p \right)=-x-\frac{1}{2}p^2,\\
\end{array}
\end{equation}
transforms the canonical form $\bar{u}'''=0$ to the canonical form (\ref{n2}).
\end{example}
\begin{example}[A Schwarzian equation associated with free-particle dynamics]
Consider the third-order nonlinear equation
\begin{equation}\label{n1}
u'''=\frac{3}{2}\frac{u''^2}{u'},
\quad u'\neq 0.
\end{equation}
The Schwarzian derivative of $u$ with respect to $x$ is defined by
$
\{u,x\}
=
\frac{u'''}{u'}
-\frac{3}{2}
\left(\frac{u''}{u'}\right)^2.
$
Consequently, equation~\eqref{n1} is equivalent to the vanishing
Schwarzian condition
$
\{u,x\}=0.
$
This condition arises naturally in projective differential geometry
and in the projective description of free-particle dynamics. More
generally, the Schwarzian equation
$
\{u,x\}=2\omega^2
$
is associated with a projective reparametrization relating the free
particle to the harmonic oscillator. Therefore, equation~\eqref{n1}
corresponds to the zero-frequency case $\omega=0$; see
\cite{Ovsienko}.

If $u$ is interpreted as a position variable and $x$ as time, then
$u'$, $u''$, and $u'''$ represent the velocity, acceleration, and
jerk, respectively. Equation~\eqref{n1} may therefore be interpreted
as a nonlinear jerk equation in which the jerk is determined by a
nonlinear velocity--acceleration relation:
$
j=\frac{3}{2}\frac{A^2}{V},
\quad
V=u',\quad A=u'',\quad j=u'''.
$
Equivalently,
$
\frac{dA}{dx}
=
\frac{3}{2}\frac{A^2}{V}.
$
This interpretation should be understood as a geometric jerk
description of the Schwarzian equation, rather than as a standard
constitutive model of a particular mechanical device.

For
$
f(x,u,p,q)=\frac{3q^2}{2p},
$
one can verify that the hypotheses of Theorem~\ref{Th} are satisfied.
Hence, equation~\eqref{n1} admits a ten-dimensional contact symmetry Lie algebra and is contact-equivalent to the canonical linear equation
$
\bar{u}'''=0.
$
For equation~\eqref{n1}, the fundamental quantities are
$
I_1=-\frac{3q}{p},
\quad
I_2=0.
$
The construction of the corresponding contact transformation proceeds
as follows.

\begin{description}
\item[Step 1.]
A solution of system~\eqref{4.21} is
$
v=1,\quad w=0.
$

\item[Step 2.]
A solution of system~\eqref{3.4} is
$
a_3=\frac{\sqrt{p}}{q}.
$

\item[Step 3.]
The auxiliary functions are $a_1=1$ and $a_2=0$.

\item[Step 4.]
The auxiliary functions are $a_7=0$ and $a_8=\frac{1}{\sqrt{p}}.$ 

\item[Step 5.]
The remaining auxiliary function is
$
r_2=0.
$

\item[Step 6.]
The auxiliary functions satisfy the compatibility equation (\ref{4.21*}) and the remaining equations of the
original auxiliary system (\ref{2.31}).

\item[Step 7.]
A solution of system~\eqref{4.22} is
$
\eta=-\frac{x}{2}-\frac{p}{q}.
$

\item[Step 8.]
A solution of system~\eqref{4.23} is
$
\chi=-x\sqrt{p}.
$

\item[Step 9.]
A solution of system~\eqref{4.24} is
$
\varphi=2\sqrt{p}.
$

\item[Step 10.]
A solution of system~\eqref{4.25} is
$
\psi=u-px.
$
\end{description}

Thus, the proper contact transformation
\begin{equation}\label{4.15}
\bar{x}=\varphi(x,u,p)=2\sqrt{p},
\qquad
\bar{u}=\psi(x,u,p)=u-px,
\qquad
\bar{p}=\chi(x,u,p)=-x\sqrt{p},
\end{equation}
transforms the nonlinear equation~\eqref{n1} into the canonical linear
equation
$
\bar{u}'''=0.
$

Consequently, transformation~\eqref{4.15} provides an explicit procedure
for obtaining the general solution of the nonlinear Schwarzian
equation~\eqref{n1} from the general solution of the canonical linear
equation $\bar{u}'''=0$. 

The significance of these two examples lies in demonstrating that although the equations are highly nonlinear, the invariant characterization established in this paper proves that they belong to the maximal contact symmetry class. Moreover, the explicit construction of the contact transformations obtained here reduces these nonlinear equations to the canonical equation $\bar u'''=0$, thereby providing a constructive linearization procedure and an effective method for generating their exact solutions.
\end{example}

\begin{example}
Consider the nonlinear third-order ODE
\begin{equation}\label{4.20}
u'''=\dfrac{3u'u''^{2}}{1+u'^{2}}.
\end{equation}
It is shown in \cite{Omar2025}  that the point transformation $\bar{x}=u+ix$, $\bar{u}=x+iu$ maps the canonical equation 
$
\bar{u}'''=\frac{3}{2}\frac{\bar{u}''^{2}}{\bar{u}'}$
into the third-order ODE $u'''=\frac{3u'u''^{2}}{1+u'^{2}}$. Hence, from the previous example, it follows that the contact transformation $$\bar{x}=2\sqrt{\frac{1+u'}{1+iu'}}, \quad \bar{u}=-2\frac{x u'-u}{1+iu'}, \bar{p}
=
-\sqrt{\frac{i+p}{1+ip}}\,(x+iu)$$
transforms the trivial equation $\bar{u}'''=0$ into the ODE (\ref{4.20}).
\end{example}
\begin{remark}
Any third-order ODE that is equivalent  under a point transformation to either of the canonical forms
$
u'''=u''^{3}$ 
or
$u'''=\frac{3}{2}\frac{u''^{2}}{u'}
$
admits a contact transformation that maps it to the canonical equation
$
u'''=0.
$
\end{remark}
\begin{example}
Consider the nonlinear ODE
\begin{equation}\label{n211}
u'''=\frac{3u''(4x^3u''-2x^2u'u''-xu'^2+3xuu''+uu')}{(2xu'-u)(xu'-u)}.\\
\end{equation}
One can verify that the function
$f(x,u,p,q)=\frac{3q(4x^3q^2-2x^2pq-xp^2+3xuq+up)}{(2xp-u)(xp-u)}$ satisfies the conditions of Theorem \ref{Th}.
Consequently, this equation admits a ten-dimensional contact symmetry Lie algebra. Moreover, it is equivalent to the canonical form
$\bar{u}'''=0$.
Using $I_1,~I_2$ for $u'''=\frac{3u''(4x^3u''-2x^2u'u''-xu'^2+3xuu''+uu')}{(2xu'-u)(xu'-u)}$, we outline the steps.
\begin{description}
\item[Step 1.]  A solution of system (\ref{4.21}) is $v=\sqrt{\frac{2xp-u}{xp-u}}, w=\frac{xp^2-xuq-up}{2(xp-u)^{3/2}\sqrt{2xp-u}}$. 
 \item[Step 2.]  A solution of system (\ref{3.4}) is $a_3=\frac{xp-u}{q(2xp-u)^2}.$
\item[Step 3.]  The auxiliary functions are $a_1=\frac{xp-u}{2xp-u}$and $a_2=-\frac{xp^2-xuq-up}{q(2xp-u)^3}.$
\item[Step 4.] The auxiliary functions are $a_7=-p$ and  $a_8=2xp-u$.
\item[Step 5.] The auxiliary function is $r_2=
\frac{
-12x^4u q^3
-3x^2u(3u-2xp)q^2
+u(xp-u)(u+2xp)q
-p^2(xp-u)^2
}{
q(2xp-u)(xp-u)^3
}.$
\item[Step 6.] The auxiliary functions satisfy the compatibility equation (\ref{4.21*}) and the remaining equations of the
original auxiliary system (\ref{2.31}).
\item[Step 7.] A solution of system (\ref{4.22}) is $\eta=\frac{-2x^2q+xp-u}{q(2xp-u)^3}$.
\item[Step 8.] A solution of  system (\ref{4.23}) is $\chi=-\frac{x}{2xp-u}$.
\item[Step 9.] A solution of  system (\ref{4.24}) is $\varphi =p(px-u)$.
\item[Step 10.] A solution of  system (\ref{4.25}) is $\psi=u-px$.
\end{description}
Finally, one can verify that the \textbf{proper contact transformation} 
\begin{equation}\label{n4}
\begin{array}{lll}
\varphi \left( x,u,p \right)=p(px-u),&\psi \left( x,u,p \right)=u-px,&\chi\left( x,u,p \right)=-\frac{x}{2xp-u},\\
\end{array}
\end{equation}
transforms the canonical form $\bar{u}'''=0$ to the nonlinear ODE (\ref{n211}).
\end{example}
\begin{remark}
The functions $f(x,u,p,q)$ in the previous examples do not satisfy the conditions of Theorem 2.1 in \cite{DweikTahirFazal2}. Therefore, equations (\ref{n2}), (\ref{n1}), (\ref{4.20}) and (\ref{n211}) are not equivalent to the canonical form $\bar{u}'''=0$ 
 under a point transformation.
\end{remark}

\section{Conclusion}
We have shown in this work how the Cartan equivalence method can be used to find an invariant characterization of scalar third-order ODEs
that admit the ten-dimensional contact symmetry algebra. Importantly this approach provides auxiliary functions which can be effectively utilized
to construct the contact transformation in order to find the reduction to the simplest third-order ODE. We have demonstrated the utility of the method by means of several significant examples presented in a constructive manner.
\subsection*{Acknowledgements}
The authors acknowledge Birzeit University for its support and facilities.\\\\
Conflict of Interest: The authors report no conflicts of interest.


\begin{thebibliography}{99}
\bibitem{Lie1}  Lie, S., Over en classe geometriske Transformationer, Doctoral Thesis, University of Christiana, 1871.
\bibitem{Lie2}  Lie, S., Begr\"undung einer Invariantentheorie der Ber\"hrungstransformationen, Mathematische Annalen 8,
               1874, 215-288.
\bibitem{Lie3}  Lie, S. and Engel, F., Theorie der Transformationsgruppen, B. G. Teubner, Leipzig, Vol. 1, 1888.
\bibitem{Lie4}  Lie, S. and Engel, F., Theorie der Transformationsgruppen, B. G. Teubner, Leipzig, Vol. 2, 1890.
\bibitem{Lie5} Lie, S. and Engel, F., Theorie der Transformationsgruppen, B. G. Teubner, Leipzig, Vol. 3, 1893.
\bibitem{Wafo} Wafo Soh C., Mahomed F. M. and Qu C, Contact Symmetry Algebras of Scalar Ordinary Differential
                Equations, Nonlinear Dynamics, 28, 213-230, 2002.
\bibitem{Lie6} Lie, S. and Scheffers, G., Vorlesugen \"uber Differentialgleichungen mit bekanten infinitesimalen Transformationen,
                B. G. Teubner, Leipzig, 1891.
\bibitem{Svi1} Svishchevskii, S. R., Lie-B\"acklund symmetries of linear ODEs and invariant linear spaces, in Modern
                Group Analysis, G. N. Yakovenko (ed.), Institute for Mathematical Modelling, Russian Academy of
                Sciences, Moscow, 1993, pp. 3-24.
\bibitem{Svi2}  Svishchevskii, S. R., Lie-B\"acklund symmetries of linear ODEs and generalized separation of variables in
                nonlinear equations, Physics Letters A 199, 1995, 344-348.
\bibitem{Ibr}   Ibragimov, N. H., Khalique, C. M. and Mahomed, F. M., All linear ordinary differential equations admitting
                contact symmetries, in Proceedings of the International Conference at the Sophus Lie Centre, N. H.
               Ibragimov, K. R. Naqvi, and E. Straume (eds.), Mars Publishers, Symmetri Foundation, Trondheim, 1997,
               pp. 155-159.
\bibitem{Yum1}  Yumaguzhin, V. A., Contact classification of 3rd-order linear ODEs, The Diffeity Institute Preprint Series,
                http://ecfor.rssi.ru/diffeity, 1997.
\bibitem{car} Cartan, E.,  {\it Bull. Soc. Math. France} {\bf52}, (1924), 205-41.
\bibitem{Gardner1989}Gardner, R. B.,  The Method of Equivalence and Its Applications, Philadelphia, SIAM, 1989.
\bibitem{Olver1995} Olver, P. J., Equivalence, Invariants and Symmetry, Cambridge University Press, Cambridge, 1995.
\bibitem{Yum}Yumaguzhin, Valeriy A., Classification of 3rd order linear ODE up to equivalence, Differential Geometry and its Applications 6. 4 (1996): 343-350.
\bibitem{mah4} Mahomed, F. M. and Leach P. G. L., Symmetry Lie Algebras of $n$th Order Ordinary Differential Equations. {\it J Math Anal Applic}  {\bf151}, (1990), 80.
\bibitem{che} Chern, S. S., The geometry of the differential equation $y'''=F(x,y,y',y'')$, Sci. Rep. Nat. Tsing Hua Univ. 4 (1940), 97-111.
\bibitem{neu} Neut, S. and Petitot, M., La g\'eom\'etrie de l'\'equation $y'''=f(x,y,y',y'')$ {\it C.R. Acad. Sci. Paris S\'er I} {\bf 335}, (2002), 515-518.
\bibitem{ibr1}  Ibragimov, N. H. and Meleshko, S. V.,  Linearization of third-order ordinary differential equations by point and contact transformations, {\it J. Math. Anal. Appl.} {\bf308}, (2005), 266-289.
\bibitem{Wunschmann} W$\ddot{\textrm{u}}$nschmann, K., $\ddot{\textrm{U}}$ber Beruhrungsbedingungen bei Differentialgleichchungen, Enzyklop$\ddot{\textrm{a}}$die der Math. Wiss. {\bf3}, (1905), 490-492.
\bibitem{DweikTahirFazal1} Al-Dweik,  Ahmad Y.,  Mahomed, F. M. and  Mustafa, M. T., Invariant characterization of third-order ODEs $u'''=f(x,u,u',u'')$ that admit a five-dimensional point symmetry Lie algebra, Communications in Nonlinear Science and
Numerical Simulation, Vol. 67, Feb, 627-636.
\bibitem{DweikTahirFazal2} Al-Dweik, Ahmad Y.,  Mustafa, M. T. and Mahomed, F. M.,  Invariant characterization of scalar third-order ODEs that admit the maximal point symmetry Lie algebra,    Mathematical Methods in the Applied Sciences, Vol. 41 Issue 12, Apr, 4714-4723, https://doi.org/10.1002/mma.4923.
\bibitem{DweikTahirFazal3} Al-Dweik, Ahmad Y., Mustafa,  M. T., Mahomed,  F. M. and Alassar, R. S.,  Linearization of third-order ordinary differential equations $u'''=f(x,u,u',u'')$ via point transformations, Mathematical Methods in the Applied Sci-
ences, Vol. 41 Issue 16, Aug, 6017-7098, https://doi.org/10.1002/mma.5208.
\bibitem{Omar2025}Abuloha, Omar A., Marwan Aloqeili, Ahmad Y. Al-Dweik, and F.
M. Mahomed. ”Equivalence Problem for Non-Linearizable Third-Order
ODEs with Four-Dimensional Lie Symmetry Subalgebras under Point
Transformations.” arXiv preprint arXiv:2602.13317 (2026).
\bibitem{Ovsienko}V. Ovsienko and S. Tabachnikov, Projective differential geometry old and new,
Cambridge University Press, Cambridge, UK, ISBN 9780521831864 (2005),
doi:10.1017/CBO9780511543142.
\end{thebibliography}
\end{document}